
\input amstex.tex
\documentstyle{amsppt}
\magnification1200
\hsize=12.5cm
\vsize=18cm
\hoffset=1cm
\voffset=2cm

\def\DJ{\leavevmode\setbox0=\hbox{D}\kern0pt\rlap
{\kern.04em\raise.188\ht0\hbox{-}}D}
\def\dj{\leavevmode
 \setbox0=\hbox{d}\kern0pt\rlap{\kern.215em\raise.46\ht0\hbox{-}}d}

\def\txt#1{{\textstyle{#1}}}
\baselineskip=13pt
\def\hf{{\textstyle{1\over2}}}
\def\b{\beta}
\def\d{{\,\roman d}}
\def\e{\varepsilon}
\def\f{\varphi}

\def\k{\kappa}
\def\s{\sigma}
\def\t{\theta}
\def\={\;=\;}

\def\D{\Delta}
\def\no{\noindent}
 
\def\z{\zeta}

 \def\t{\theta}
\def\hf{{\textstyle{1\over2}}}
\def\txt#1{{\textstyle{#1}}}
\def\f{\varphi}

\font\tenmsb=msbm10
\font\sevenmsb=msbm7
\font\fivemsb=msbm5
\newfam\msbfam
\textfont\msbfam=\tenmsb
\scriptfont\msbfam=\sevenmsb
\scriptscriptfont\msbfam=\fivemsb
\def\Bbb#1{{\fam\msbfam #1}}

\def \NN {\Bbb N}

\def \RR {\Bbb R}
\def \ZZ {\Bbb Z}

\font\ff=cmr8
\def\txt#1{{\textstyle{#1}}}
\baselineskip=13pt

\font\teneufm=eufm10
\font\seveneufm=eufm7
\font\fiveeufm=eufm5
\newfam\eufmfam
\textfont\eufmfam=\teneufm
\scriptfont\eufmfam=\seveneufm
\scriptscriptfont\eufmfam=\fiveeufm
\def\mathfrak#1{{\fam\eufmfam\relax#1}}

\font\tenmsb=msbm10
\font\sevenmsb=msbm7
\font\fivemsb=msbm5
\newfam\msbfam
     \textfont\msbfam=\tenmsb
      \scriptfont\msbfam=\sevenmsb
      \scriptscriptfont\msbfam=\fivemsb
\def\Bbb#1{{\fam\msbfam #1}}

\def \NN {\Bbb N}

\def \RR {\Bbb R}
\def \ZZ {\Bbb Z}

  \def\rightheadline{{\hfil{\ff
On the fourth moment in the Rankin-Selberg problem }\hfil\tenrm\folio}}

  \def\leftheadline{{\tenrm\folio\hfil{\ff
   Aleksandar Ivi\'c }\hfil}}
  \def\emptyheadline{\hfil}
  \headline{\ifnum\pageno=1 \emptyheadline\else
  \ifodd\pageno \rightheadline \else \leftheadline\fi\fi}

\font\ff=cmr8
\font\teneufm=eufm10
\font\seveneufm=eufm7
\font\fiveeufm=eufm5
\newfam\eufmfam
\textfont\eufmfam=\teneufm
\scriptfont\eufmfam=\seveneufm
\scriptscriptfont\eufmfam=\fiveeufm
\def\mathfrak#1{{\fam\eufmfam\relax#1}}

\font\tenmsb=msbm10
\font\sevenmsb=msbm7
\font\fivemsb=msbm5
\newfam\msbfam
\textfont\msbfam=\tenmsb
\scriptfont\msbfam=\sevenmsb
\scriptscriptfont\msbfam=\fivemsb
\def\Bbb#1{{\fam\msbfam #1}}

\def \NN {\Bbb N}

\def \RR {\Bbb R}
\def \ZZ {\Bbb Z}

\def\D{\Delta}

\def\b{\beta} \def\e{\varepsilon}
\def\no{\noindent} \def\d{\,{\roman d}}
\topmatter
\title
On the fourth moment in the Rankin-Selberg problem
\endtitle
\author
Aleksandar Ivi\'c
\endauthor
\address
Katedra Matematike RGF-a, Universitet u Beogradu,  \DJ u\v sina 7,
11000 Beograd, Serbia
\endaddress
\subjclass 11 N 37, 11 M 06, 44 A 15, 26 A 12
\endsubjclass
\email {\tt  ivic\@rgf.bg.ac.yu, aivic\@matf.bg.ac.yu}
\endemail
\abstract
If
$$
\D(x) \;:=\; \sum_{n\le x}c_n - Cx
$$
denotes the error term in the classical Rankin-Selberg
problem, then it is proved that
$$
\int_0^X \D^4(x)\d x \ll_\e X^{3+\e},\quad
\int_0^X \D_1^4(x)\d x \ll_\e X^{11/2+\e},
$$
where $\D_1(x) = \int_0^x\D(u)\d u$. The latter bound is,
up to `$\e$', best possible.

\endabstract
\endtopmatter
\document

\head 1. Introduction and statement of results
\endhead
The classical Rankin-Selberg problem consists
of the estimation of the error term function
$$
\D(x) \;:=\; \sum_{n\le x}c_n - Cx,\leqno(1.1)
$$
where the notation is as follows. Let $\varphi(z)$ be a
holomorphic cusp form of weight $\kappa$ with respect to the full
modular group $SL(2,\ZZ)$, and denote by $a(n)$ the $n$-th Fourier
coefficient of $\varphi(z)$. We suppose that $\varphi(z)$ is a
normalized eigenfunction for the Hecke operators $T(n)$, that is,
$  a(1)=1  $ and $  T(n)\varphi=a(n)\varphi $ for every $n \in
\NN$. The classical example is $a(n) = \tau(n)\;(\k=12)$, the
Ramanujan function defined by
$$
\sum_{n=1}^\infty \tau(n)x^n \=
x{\left\{(1-x)(1-x^2)(1-x^3)\cdots\right\}}^{24}\qquad(\,|x| < 1).
$$
The constant $C (>0)$ in (1) may be written down explicitly
(see e.g., [7]), and $c_n$ is the convolution function
defined by
$$
c_{n}=n^{1-\kappa}\sum_{m^2 \mid n}m^{2(\kappa-1)}
\left|a\Bigl({n\over m^2}\Bigr)\right|^2.
$$
The classical Rankin-Selberg  bound of 1939 is
$$
\D(x) = O(x^{3/5}),\leqno(1.2)
$$
hitherto unimproved. In their works, done independently, R.A.
Rankin [10] derives (1.2) from a general result of E. Landau [9],
while A. Selberg [12] states the result with no proof. Although
it seems very difficult at present to improve the bound in (1.2),
recently there have been some results on mean square estimates in
the Rankin-Selberg problem (see the author's works [5],  [6]).
Namely, let as usual $\mu(\s)$ denote the Lindel\"of function
$$
\mu(\s) := \limsup_{t\to\infty}\,{\log|\z(\s+it)|\over\log
t}\qquad(\s\in\RR).
$$
Then we have
$$
\int_0^X \D^2(x)\d x \ll_\e X^{1+2\b+\e},\quad\b \= {2\over 5-2\mu(\hf)}.\leqno(1.3)
$$
Here and later $\e$ denotes positive
constants which may be arbitrarily small, but are not necessarily the same at
each occurrence, while $\ll_\e$ means that the $\ll$--constant depends on $\e$.
Note that with the sharpest known  result (see M.N.
Huxley [2]) $\mu(\hf) \le 32/205$ we obtain $\b = 410/961 =
0.426638917\ldots\;$. The limit of (1.3) is the value $\b = 2/5$
if the Lindel\"of hypothesis (that $\mu(\hf) =0$) is true.

\medskip
We propose to contribute
here to the subject of mean value results  for $\D(x)$
by proving (unconditionally) the following results.

\bigskip
THEOREM 1. {\it For any given $\e>0$ we have}
$$
\int_0^X \D^4(x)\d x \ll_\e X^{3+\e}.\leqno(1.4)
$$

\medskip\no
Note that (1.4) follows from (1.3) only if the
Lindel\"of hypothesis  $\mu(\hf) =0$ is true.
\medskip
{\bf Corollary}. {\it For any given $\e>0$ we have}
$$
\D(x) \ll_\e x^{3/5+\e}.\leqno(1.5)
$$

\medskip
\no
The bound in (1.5) is only by an `$\e$'--factor weaker than the
strongest known bound (1.2). To obtain (1.5) from (1.4) note that
we have (see Lemma 1 below)
$$
\D(X) = {1\over2H}\int_{X-H}^{X+H}\D(x)\d x + O(H)\qquad(X^\e\le H \le
\hf X). \leqno(1.6)
$$
It follows from (1.6) by H\"older's inequality for integrals that
$$
\D^4(X) \ll H^{-1}\int_{X-H}^{X+H}\D^4(x)\d x  + H^4 \ll_\e H^{-1}X^{3+\e} + H^4
$$
by (1.4), and (1.2) follows with $H = X^{3/5}$.
Note that if (1.4) holds with the exponent $\t$ on the right-hand side of (1.4),
then the above argument gives
$$
\D(x) \ll x^{\t/5},
$$
and the best possible exponent $\t$ must satisfy $\t \ge 15/8$ since
$\D(x) = \Omega_\pm(x^{3/8})$ (see the author's work [4]).

\bigskip
THEOREM 2. {\it If $\D_1(x) = \int_0^x\D(u)\d u$, then for any given $\e>0$ we have}
$$
\int_0^X \D_1^4(x)\d x \ll_\e X^{11/2+\e}.\leqno(1.6)
$$

\medskip
Note that it was proved in [7] that
$$
\int_0^X\D_1^2(x) \d x = {2\over13}(2\pi)^{-4}\left(\sum_{n=1}^\infty
c_n^2n^{-7/4}\right)X^{13/4} + O_\e(X^{3+\e}),\leqno(1.7)
$$
so that from (1.7) and the Cauchy-Schwarz inequality for integrals we obtain that
$$
\int_0^X \D_1^4(x)\d x \;\gg\; X^{11/2}.
$$
This shows that, up to `$\e$', the bound in (1.6) is best possible.

\bigskip
\head
2. The necessary lemmas
\endhead
In this section we shall state the lemmas necessary for the proof of our theorems.

\medskip
LEMMA 1. {\it For $X^\e\le H \le \hf X$ we have}
$$
\D(X) = {1\over2H}\int_{X-H}^{X+H}\D(x)\d x + O(H).\leqno(2.1)
$$

\medskip
{\bf Proof.}
$$
\eqalign{
\D(X) - {1\over2H}\int_{X-H}^{X+H}\D(x)\d x &
= {1\over2H}\int_{X-H}^{X+H}(\D(X) - \D(x))\d x\cr&
= {1\over2H}\int_{X-H}^{X+H}\left(\sum_{x<n\le X}c_n + C(X-x)\right)\d x\cr&
\ll \sum_{X-H\le n\le X+H}c_n + H\ll H,\cr}
$$
where in the last step a well-known result of P. Shiu on multiplicative functions
[14] in short intervals was used (see also Lemma 4 of [7]).

\medskip
The next two lemmas are the explicit, truncated formula
of the Vorono{\"\i} type for  for $\D(x)$ and $\D_1(x)$, respectively.

\medskip
LEMMA 2. {\it For $1 \ll K_0 \ll x$ a parameter we have}
$$
\D(x) = {x^{3/8}\over2\pi}\sum_{k\le K_0}c_kk^{-5/8}\sin
\left(8\pi(kx)^{1/4}+{\txt{3\pi\over4}}\right) +
O_\e\Bigl(x^{3/4+\e}K_0^{-1/4}\Bigr), \leqno(2.2)
$$

\medskip\no
Choosing $K_0 = x^{3/5}$ and  estimating
trivially the sum in (2.2) we obtain again the bound $\D(x) \ll_\e
x^{3/5+\e}$. The formula
(2.2) was proved (see [7, Lemma 2] with $\rho=0$)
by Ivi\'c, Matsumoto and Tanigawa, where
the general (Riesz) sum $\sum_{n\le x}(x-n)^\rho c_n$ for fixed
$\rho \ge 0$ is investigated, and a proof of (1.2) is given. A
similar formula holds in general for Dirichlet series of degree
four in the Selberg class (see e.g., the survey article [8]
of Kaczorowski-Perelli on functions from the Selberg class $\Cal S$,
and Selberg's original paper [13]). Thus the generalization of our
results will hold for error terms associated to
suitable Dirichlet series in $\Cal S$. However,
in some cases the special structure of the problem at hand allows
for sharper results. For example, consider the estimation of $\D_4(x)$,
the error term in the asymptotic formula for the summatory function of
the divisor function $d_4(n) = \sum_{abcd=n;a,b,c,d\in\NN}1$.
The generating function in this case is $\z^4(s)$, and we have
(see e.g., [3, Chapter 13])
$$
\int_0^X\D^2_4(x)\d x \;\ll_\e\; X^{7/4+\e},\leqno(2.3)
$$
which is (up to `$\e$') best possible in view of
$\D_4(x) = \Omega_\pm(x^{3/8})$ (see [4]).
Since one has $\D_4(x) \ll_\e x^{1/2+\e}$, it trivially follows from (2.3)
that
$$
\int_0^X\D^4_4(x)\d x \;\ll_\e\; X^{11/4+\e},\leqno(2.4)
$$
and the exponent in (2.4) is better than the exponent in (1.4).

\medskip
LEMMA 3. {\it For $1 \ll K_0 \ll x^2$ a parameter we have}
$$
\D_1(x) = {x^{9/8}\over(2\pi)^2}\sum_{k\le K_0}c_kk^{-7/8}\sin
\left(8\pi(kx)^{1/4}+{\txt{\pi\over4}}\right) + O_\e\Bigl(x^{1+\e}+
x^{3/2+\e}K_0^{-1/2}\Bigr). \leqno(2.5)
$$

\medskip\no Lemma 3 is the case $\rho = 1$ of [7, Lemma 2].

\medskip

The chief ingredient in the proof of Theorem 2 is the
new result of O. Robert--P. Sargos [11], which is the following

\medskip
LEMMA 4. {\it Let $k\ge 2$ be a fixed
integer and $\delta > 0$ be given.
Then the number of integers $n_1,n_2,n_3,n_4$ such that
$N < n_1,n_2,n_3,n_4 \le 2N$ and}
$$
|n_1^{1/k} + n_2^{1/k} - n_3^{1/k} - n_4^{1/k}| < \delta N^{1/k}
$$
{\it is, for any given $\e>0$,}
$$
\ll_\e N^\e(N^4\delta + N^2).\leqno(2.6)
$$

\head
3. Proof of Theorem 1
\endhead
\medskip
To prove (1.4), it is sufficient to prove the bound in question
for the integral over the interval $[X,\,2X]$. Let, for
$$
X^{1/2} \;\le\; V \ll X^{3/5},\leqno(3.1)
$$
$$
I(V,X) \;:=\; \int_{X,V\le|\D(x)|<2V}^{2X}\D^4(x)\d x,\leqno(3.2)
$$
where the upper bound in (3.1) holds because of (1.2). Clearly we have
$$
\int_X^{2X}\D^4(x)\d x \ll X^3 + \log X\max_{X^{1/2}\le V\ll X^{3/5}}
I(V,X).
\leqno(3.3)
$$
Hence the problem is reduced to the estimation of $I(V,X)$ in (3.2).
Thus we fix a value $V = C2^{-j}X^{3/5}\,(j = 0,1,2,\ldots,\,C>0)$
and split the
interval $[X,\,2X]$  into subintervals of length $H\,(X^{1/2}
\le H \le X^{1-\e})$, where the last of these intervals may be shorter.
Suppose there are $R = R(V)$ of these subintervals which contain a point
$x$ for which $V\le|\D(x)|<2V$ holds. Further suppose that $x_r$
is the point in the $r$-th of these intervals where the largest value of
$|\D(x)|$ is attained. To obtain the spacing condition
$$
|x_r - x_s| \;\ge \;H\quad(r\ne s;\; r,s = 1,\ldots,R)\leqno(3.4)
$$
we consider separately the points with even and odd indices and,
with a slight abuse of notation, each of these two systems of points
is again denoted by $\{x_r\}_{r=1}^R$. Therefore we have
$|\D(x_r)| \ge V\,(r = 1,2\ldots\,)$, and observe that
$$
{x^{3/8}\over2\pi}\sum_{k\le \delta X^{1/3}}c_kk^{-5/8}\sin
\left(8\pi(kx)^{1/4}+{\txt{3\pi\over4}}\right) \ll \delta^{3/8}X^{1/2}
$$
for $x = x_r\in\, [X,\,2X], \,\delta>0$. Thus for $\delta$ small
enough it follows from (2.1) and (2.2)(changing $K_0$ to $X^3H^{-4}$ and
recalling that ${\roman e}(z) = {\roman e}^{2\pi iz}$), that for
$r = 1,\ldots\,,R$
$$
V \ll {X^{3/8}\over H}\int_{x_r-H/3}^{x_r+H/3}
\Bigl|\sum_{\delta X^{1/3}\le k\le X^3H^{-4}}c_kk^{-5/8}{\roman e}
(4\pi(kx)^{1/4})\Bigr|\d x + HX^\e,\leqno(3.5)
$$
where all the intervals $[x_r-H/3,\,x_r+H/3]$ are disjoint in view of (3.4).
We take in (3.5)
$$
H \;=\; VX^{-2\e},
$$
square, use the Cauchy-Schwarz inequality for integrals and sum the
resulting expressions. We obtain
$$
R \ll_\e \max_{\delta X^{1/3}\le K\le X^{3+8\e}V^{-4}}
X^{3/4+2\e}V^{-3}\int_{X/2}^{5X/2}
\Bigl|\sum_{K<k\le2K}c_kk^{-5/8}{\roman e}
(4\pi(kx)^{1/4})\Bigr|^2\d x.
$$
To evaluate the integral on the right-hand side we square out the
sum and use the first derivative test (i.e., Lemma 3.1 of [1] or Lemma
2.1 of [3]). It follows that, since $c_n \ll_\e n^\e$,
$$
\eqalign{&
\int_{X/2}^{5X/2}
\Bigl|\sum_{K<k\le2K}c_kk^{-5/8}{\roman e}
(4\pi(kx)^{1/4})\Bigr|^2\d x \cr&
 = \sum_{K <k_1,k_2\le2K}c_{k_1}c_{k_2}(k_1k_2)^{-5/8}\int_{X/2}^{5X/2}
{\roman e}\Bigl(4\pi x^{1/4}(k_1^{1/4}-k_2^{1/4})\Bigr)\d x\cr&
\ll_\e X^{1+\e}K^{-1/4} + K^{-5/4}\sum_{K <k_1\ne k_2\le2K} {X^{3/4+\e}\over
|k_1^{1/4}-k_2^{1/4}|}\cr&
\ll_\e X^{1+\e}K^{-1/4} + K^{-1/2}X^{3/4+\e}\sum_{K <k_1\le2K}
\sum_{K <k_2\le2K,k_2\ne k_1}{1\over|k_2-k_1|} \cr&
\ll_\e X^{1+\e}K^{-1/4} + X^{3/4+\e}K^{1/2} \ll_\e X^{3/4+\e}K^{1/2}\cr}
$$
for $K \gg X^{1/3}$, which is the case in (3.5). Consequently
$$
R \ll_\e \max_{\delta X^{1/3}\le K\le X^{3+8\e}V^{-4}}
X^{3/4+2\e}V^{-3}X^{3/4+\e}K^{1/2} \ll_\e X^{3+\e}V^{-5}.\leqno(3.6)
$$
From (3.6) we obtain that
$$
I(V,X) \ll_\e X^\e V\cdot V^4R \ll_\e X^{3+\e},
$$
and (3.3) gives
$$
\int_X^{2X}\D^4(x)\d x \;\ll_\e\;X^{3+\e}.
$$
Theorem  1 follows if we replace $X$ by $X2^{-j}$ in the above bound
and sum over $j = 1,2,\ldots\;$. An alternative proof of (1.4) may be obtained
by going through the appropriate modification of the
proof of Theorem 13.8 in [3], taking $k=4$ and $R=R_0$,
$T = T_0$ in (13.65). The choice $(\kappa,\lambda) = (\hf,\hf)$ will provide
again the bound in (3.6) for the relevant range.

\bigskip
\head
4. Proof of Theorem 2
\endhead
\bigskip
To prove Theorem 2 we use (2.5) of Lemma 3 with $K_0 = X$ for $X/2\le x \le 5X/2$.
This gives
$$\eqalign{
\int_X^{2X}\D_1^4(x)\d x &\ll_\e X^{9/2}\log X\max_{K\ll X}\int_{X/2}^{5X/2}\f(x)
\Bigl|\sum_{K<k\le2K}c_kk^{-7/8}{\roman e}(4(xk)^{1/4})\Bigr|^4\d x \cr&+ X^{5+\e},
\cr}\leqno(4.1)
$$
where $\f(x)\;(\ge0)$ is a smooth function supported in $[X/2,\,5X/2]$ such that
$\f(x) = 1$ for $X\le x \le 2X$ and $\f^{(r)}(x) \ll_r X^{-r}\;
(r = 0,1,\ldots\,)$.
If we set
$$
\D := k^{1/4} + \ell^{1/4} - m^{1/4} - n^{1/4}\qquad(k,\ell,m,n\in\NN),
$$
then
$$
\int\limits_{X/2}^{5X/2}\f(x)\cdots\d x = \sum_{K<k,\ell,m,n\le2K}
c_kc_\ell c_mc_n(k\ell mn)^{-7/8}\int\limits_{X/2}^{5X/2}\f(x)
{\roman e}(4\pi\D x^{1/4})\d x.
$$
But integration by parts shows that
$$
\int\limits_{X/2}^{5X/2}\f(x){\roman e}(4\pi\D x^{1/4})\d x
= - \int\limits_{X/2}^{5X/2}{{\roman e}(4\pi\D x^{1/4})\over2\pi i\D}
\left(x^{3/4}\f(x)\right)'\d x.
$$
Thus the exponential factor remained the same, but the order of the integrand
has decreased by $\D^{-1}X^{-1/4}$, provided that $\D\ne0$. Since this will
be repeated after every integration by parts, then it follows that the contribution
of quadruples $(k,\ell,m,n)$ will be negligible if $\D \ge X^{\e-1/4}$ for any
given $\e>0$. The contribution of the quadruples satisfying $\D \le X^{\e-1/4}$
is estimated by  Lemma 4, where in (2.6) we take $k=4, \delta \asymp \D K^{-1/4}$.
The ensuing integral is estimated trivially (using $c_n \ll_\e n^\e$), and we obtain
$$\eqalign{
\int_X^{2X}\D_1^4(x)\d x &\ll_\e X^{9/2+\e}X\max_{K\ll X}K^{-7/2}
\Bigl(K^4(KX)^{-1/4} +K^2\Bigr) + X^{5+\e} \cr&
\ll_\e X^{21/4+\e}\max_{K\ll X}K^{1/4} + X^{11/2+\e} \ll_\e X^{11/2+\e}.\cr}
$$
Theorem  2 follows if we replace $X$ by $X2^{-j}$ in the above bound
and sum over $j = 1,2,\ldots\;$.

\vfill
\eject
\topskip2cm

\bigskip\bigskip
\Refs
\bigskip

\item{[1]} S.W. Graham and G. Kolesnik, Van der Corput's method
of exponential sums, LMS Lecture Note Series {\bf {126}},
Cambridge University Press, Cambridge, 1991.

\item{[2]} M.N. Huxley, Exponential sums and the Riemann zeta-function V,
Proc. London Math. Soc. (4) {\bf90}(2005), 1-41.

\item{[3]} A. Ivi\'c,  The Riemann Zeta-Function,
   John Wiley \& Sons, New York, 1985 (2nd ed. Dover, 2003).

\item{[4]} A. Ivi\'c, Large values of certain number-theoretic error terms,
    Acta Arith. {\bf 56}(1990), 135-159.

\item{[5]} A. Ivi\'c, Estimates of convolutions of certain
number-theoretic error terms, Intern. J. Math. and Math. Sciences,
Vol. 2004, No. {\bf1}(2004), 1-23.

\item{[6]} A. Ivi\'c, On some mean square estimates in the
Rankin-Selberg problem, Applicable Analysis
and Discrete Mathematics {\bf1}(2007), 1-11.

\item{[7]} A. Ivi\'c, K. Matsumoto and Y. Tanigawa, On Riesz mean
of the coefficients of the Rankin--Selberg series, Math. Proc.
Camb. Phil. Soc. {\bf127}(1999), 117-131.

\item{[8]} A. Kaczorowski and A. Perelli, The Selberg class: a
survey, in ``Number Theory in Progress, Proc. Conf. in honour
of A. Schinzel (K. Gy\"ory et al. eds)", de Gruyter,
Berlin, 1999, pp. 953-992.

\item{[9]}E.~Landau, \"{U}ber die Anzahl der
Gitterpunkte in gewissen Bereichen II, Nachr. Ges.
Wiss. G\"ottingen 1915, 209-243.

\item{[10]}R.~A.~Rankin, Contributions to the theory of Ramanujan's
   function $\tau(n)$ and similar arithmetical functions II. The order
   of the Fourier coefficients of integral modular forms,
    Proc. Cambridge Phil. Soc. {\bf 35}(1939), 357-372.

\item{[11]} O. Robert and P. Sargos, Three-dimensional
exponential sums with monomials, J. reine angew. Math. {\bf591}(2006), 1-20.

    \item{[12]}A.~Selberg,  Bemerkungen \"{u}ber eine Dirichletsche Reihe,
    die mit der Theorie der Modulformen nahe verbunden ist,
    Arch. Math. Naturvid. {\bf 43}(1940), 47-50.

    \item{[13]} A.~Selberg, Old and new conjectures and results about a class
    of Dirichlet series, in ``Proc. Amalfi Conf. Analytic Number Theory 1989
    (E. Bombieri et al. eds.)",
    University of Salerno, Salerno, 1992, pp. 367--385.

\item{[14]} P. Shiu, A Brun--Titchmarsh theorem for multiplicative functions,
J. reine angew. Math. {\bf313}(1980), 161-170.

\vskip1cm
\endRefs

\enddocument

\end